\documentclass[12pt]{article}

\usepackage{amsthm}
\usepackage{amssymb}
\usepackage{amsmath}
\usepackage{hhline}
\usepackage[mathscr]{eucal}

\ifx\pdfoutput\undefined
   \usepackage[dvips]{graphicx}
   \else
   \usepackage[pdftex]{graphicx}
   \fi
   \usepackage{epstopdf}

\theoremstyle{plain}
\newtheorem{tetel}{Theorem}%[section]

\newtheorem{lemma}[tetel]{Lemma}
\newtheorem{kov}[tetel]{Corollary}

\theoremstyle{definition}\newtheorem{Def}[tetel]{Definition}
\theoremstyle{remark}\newtheorem{megj}[tetel]{Remark}
\newtheorem{pelda}[tetel]{Example}

\newcommand*{\R}{\ensuremath{\mathbf R}}
\newcommand*{\di}{\mathrm d}

\newcommand*{\TB}{\overline{tb}}

\begin{document}

\title{Maximal Thurston--Bennequin number of +adequate links}
\author{Tam\'as K\'alm\'an\\ University of Southern California}
\maketitle

\begin{abstract}
The class of $+$adequate links contains both alternating and positive links. Generalizing results of Tanaka (for the positive case) and Ng (for the alternating case), we construct fronts of an arbitrary $+$adequate link $A$ so that the diagram has a ruling, therefore its Thurston--Bennequin number is maximal among Legendrian representatives of $A$. We derive consequences for the Kauffman polynomial and Khovanov homology of $+$adequate links.
\end{abstract}

Maximum Thurston--Bennequin number, denoted by $\TB$, is a knot invariant 
%closely related to genus, 
that has drawn a lot of recent interest. Its definition is possible because Bennequin's inequality, $tb\le 2g-1$, bounds from above the Thurston--Bennequin number of Legendrian representatives\footnote{To avoid undue repetition and to keep this note short, for the standard definitions of Legendrian knot theory we refer the reader to \cite{etn}, \cite{en}, or to any number of other publications. Let us only state that we work in $\R^3_{xyz}$ where the contact structure is the kernel of $\di z-y\di x$, so that the front projection is the $xz$--projection.} of any knot type by (essentially) the genus $g$ of the knot. Either Bennequin's inequality itself or other bounds, for example the so-called Kauffman bound on $tb$ \cite{congr}, make the extension to links possible.

%%%%%%%%%%%%%%%%%%%%%%%%%%%%%%
% Tanaka references Tabachnikov for the Kauffman bound!
%%%%%%%%%%%%%%%%%%%%%%%%%%%%%%

Recall that the Thurston--Bennequin number is computed from an oriented front diagram by subtracting the number of right cusps from the writhe of the diagram; $\TB$ is the maximum of these numbers for all fronts representing a given link type. The Kauffman bound states that $tb$, and thus $\TB$ is strictly less than the minimum $v$--degree (or $-1$ times the maximum $a$--degree) of the Kauffman polynomial\footnote{It seems to be standard to use either $v$ and $z$ or $a=v^{-1}$ and $z$ as the variables in the Kauffman polynomial; see \cite{crom}. The `Dubrovnik version' has the same degree distribution.}.

The value of $\TB$ 
%(at least implicitly, through constructions) 
is known for the following infinite classes of knots and links: positive links \cite{tanaka} (see also \cite{yok}) 2-bridge links \cite{2bridge} and more generally, alternating links \cite{khov}, negative torus knots \cite{EH1}, and Whitehead doubles with sufficiently negative framing \cite{fuchs}.

In the positive and alternating cases, the proofs of Tanaka and Ng proceed as follows: For the given knot or link, they construct a certain front diagram. Tanaka uses the Kauffman bound (i.e., establishes that it's sharp for the front) to show that its Thurston--Bennequin number is maximal. Ng uses the so-called Khovanov bound for the same purpose, but it turns out that the Kauffman bound works just as easily in his case as well. (The Kauffman bound is known to be not sharp for many other knots, such as those where the Khovanov bound actually improves it \cite{khov}, and also for negative $(p,q)$ torus knots with $p>q$ and $q$ even\footnote{What Fuchs proves about these is that they do not possess rulings; see below for the significance of this fact in the light of Rutherford's theorem. The same observation also shows that it is not possible to extend the methods of this paper to homogeneous links.} \cite{fuchs}.) Recently, Rutherford \cite{rulpoly} clarified when the Kauffman bound is sharp by verifying the Fuchs conjecture: a necessary and sufficient condition is the existence of a so-called ungraded ruling (a.k.a.\ decomposition) for at least one Legendrian representative of the link. That representative then has maximum Thurston--Bennequin number. 

\begin{figure}
   \centering
   \includegraphics[width=\linewidth]{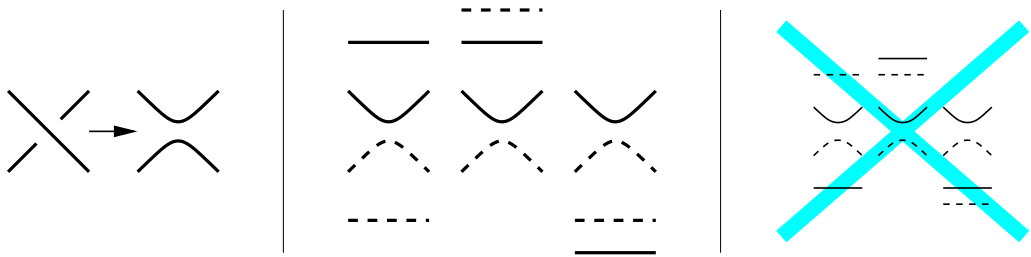}
   \caption{Allowed and disallowed configurations for switches of rulings}
   \label{fig:switch}
\end{figure}

\begin{Def}\label{def:rul} An \emph{ungraded ruling} is a partial splicing
of a front diagram where certain crossings, called \emph{switches}, are
replaced by a pair of arcs as in Figure \ref{fig:switch} so that the
diagram becomes a (not necessarily disjoint) union of standard unknot diagrams, called \emph{eyes}.
(An eye is a pair of arcs connecting the same two cusps that contain no
other cusps and that otherwise do not meet, not even at switches.) We also impose the so-called normality condition: in the vertical ($x=\text{const.}$) slice of the diagram
through each switch, the two eyes that meet at the switch follow one of
the three configurations in the middle of Figure \ref{fig:switch}.
\end{Def}

The notion of a $+$adequate link is a common generalization of positive and alternating links\footnote{Another such generalization, homogeneous links, will be alluded to in Example \ref{ex:KT}.}, introduced in \cite{adeq} (see also \cite{thist}). In this paper, we use Rutherford's theorem to unify\footnote{We could use Proposition 7 of \cite{khov} to the same effect.} Tanaka's and Ng's approaches: that is, we construct a Legendrian representative with a ruling for each $+$adequate link.

\begin{Def}
Let $D$ be a link diagram. Let $s_+(D)$ denote the diagram obtained from $D$ by splicing each crossing as on the left of Figure \ref{fig:switch}. 
This is a disjoint union of simple closed curves, called state circles. $D$ is \emph{$+$adequate} if in $s_+(D)$, the two segments replacing each crossing of $D$ belong to different circles. A link is $+$adequate if it possesses a $+$adequate diagram.
\end{Def}

An alternating diagram in which every separating crossing is positive is easily seen to be $+$adequate: the state circles are the perimeters of the black regions of a checkerboard coloring that locally looks like \includegraphics%[height=.8\lineheight]
{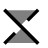} at every crossing. In particular, none of the state circles contains any other in its interior. Positive links are $+$adequate because the state circles derived from a positive diagram agree with the Seifert circles. (Here, a positive diagram is one whose algebraic and geometric crossing numbers agree. Positive links are not to be confused with the more restrictive notion of a braid-positive link, i.e.\ a link that can be obtained as the closure of a positive braid.)

\begin{figure}[b]
   \centering
   \begin{minipage}[c]{.6\textwidth}
   \centering
   \includegraphics[width=\textwidth]{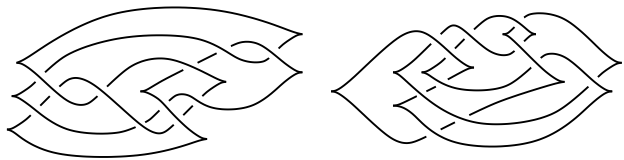}
   \end{minipage}
   \hfill
   \begin{minipage}[c]{.3\textwidth}
   \centering
   \caption{Front diagrams for the inadequate knot $11n_{95}$ and its mirror.}
   \label{fig:nemadekvat}
   \end{minipage}
\end{figure}

\begin{megj}
A knot doesn't have to be $+$adequate to have a ruling. Figure \ref{fig:nemadekvat} shows front diagrams of the inadequate knot $11n_{95}$ \cite[p.\ 230]{crom} and its mirror that realize the Kauffman bound, hence have rulings (and demonstrate that $\TB(11n_{95})=3$ and $\TB(11n^*_{95})=-12$).
%Fuchs \cite{fuchs} constructs front diagrams of `sufficiently positively twis\-ted' Whitehead doubles with a ruling. Note however that  Whitehead doubles are usually not $+$adequate.
%Is that really true??? Proof: leading coefficient in the Jones polynomial?? Or negative coefficients in Kauffman??
What is true in the converse direction is that if in a front diagram, the set of all crossings is a ruling (Ng calls these fronts with admissible $0$-resolution in \cite{khov}), then that diagram is $+$adequate.
\end{megj}

\begin{tetel}\label{front}
Let $A$ be a $+$adequate link. Then $A$ possesses a Legendrian representative with an ungraded ruling. 
\end{tetel}

In fact, any $+$adequate link diagram can be isotoped into a front diagram where the set of all crossings serves as the set of switches in an ungraded ruling. Figure \ref{fig:kinoshita} summarizes the procedure; in particular, the eyes in the ruling are exactly the state circles. (To follow the proof below, and in particular the proof of Lemma \ref{lem:vazze}, we suggest to first read pages 1647--1649 of \cite{khov}, in particular steps 1, 2, and 3.)

\begin{proof}
%the former, with its appeal to Yamada's bunching operations \cite{yam}, seems to be more straightforward. 
We may generalize either Tanaka's or Ng's construction, but the latter seems to be more straightforward. In \cite{khov}, it is proved that every alternating link projection can be deformed by a diffeomorphism of the plane into ``Mondrian position'' where the black regions are arbitrarily close to horizontal line segments, and the two arcs in a neighborhood of each crossing are arbitrarily close to vertical line segments connecting the horizontal ones. In such a position, the horizontal line segments are essentially the eyes of an ungraded ruling.

A more precise formulation of Ng's results is the following: A \emph{Mondrian diagram} is a union of a set of disjoint horizontal line segments and a set of disjoint vertical ones; each of the latter start and end on a horizontal segment, and doesn't intersect other horizontal segments. Contracting a Mondrian diagram means contracting the horizontal segments to points; Ng proves that any\footnote{Proposition 11 of Ng's paper claims this result for reduced planar graphs (i.e., ones with no separating edges) but he also remarks that the proof can be easily extended to the general case.} planar graph is the contraction of some Mondrian diagram, and then applies this to the graph whose vertices are the black regions of the alternating diagram and whose edges correspond to crossings.

In the case of a general $+$adequate diagram $D$, state circles are nested, giving rise to a partial ordering among them. An example is shown in Figure \ref{fig:kinoshita}. We may first apply Ng's construction to the diagram consisting of those circles that are not contained in others, i.e.\ ones that are maximal in the partial order (these are 1, 2, 3, and 4 in Figure \ref{fig:kinoshita}). More precisely, define a planar graph whose vertices are the maximal state circles and whose edges correspond to crossings of $D$ connecting two such. Obtain this planar graph as the contraction of a Mondrian diagram. Then, the resulting horizontal line segments (corresponding to the maximal state circles) can be slightly thickened into rectangles to make room for the remaining parts of the diagram, which are arranged again by starting from the outermost state circles. This can be applied successively toward lower and lower levels of the partial order, resulting in an `iterated Mondrian diagram.' (The upper right subdiagram of Figure \ref{fig:kinoshita} is an illustration.) Below, we explain this procedure in greater detail.
%along the branches of the obvious forest graph that has the state circles as its vertices. 

An important technical difficulty is that the order of outward and inward crossings along any state circle needs to be respected. (See Figure \ref{fig:inandout}; note that another distinguishing feature of alternating diagrams is that along their state circles, with the possible exception of a single outermost circle, all crossings are outer.) In our process, we turn each state circle into a narrow horizontal rectangle and we arrange the outward crossings first along the upper and lower sides of the rectangle. Thus, a subdivision of the inward crossings is created into upper and lower ones as well (those inward crossings found at the ends of the rectangle can be arbitrarily assigned to upper or lower). 

\begin{figure}
   \centering
   %\begin{minipage}[c]{.4\textwidth}
   %\centering
   \includegraphics[width=4in]{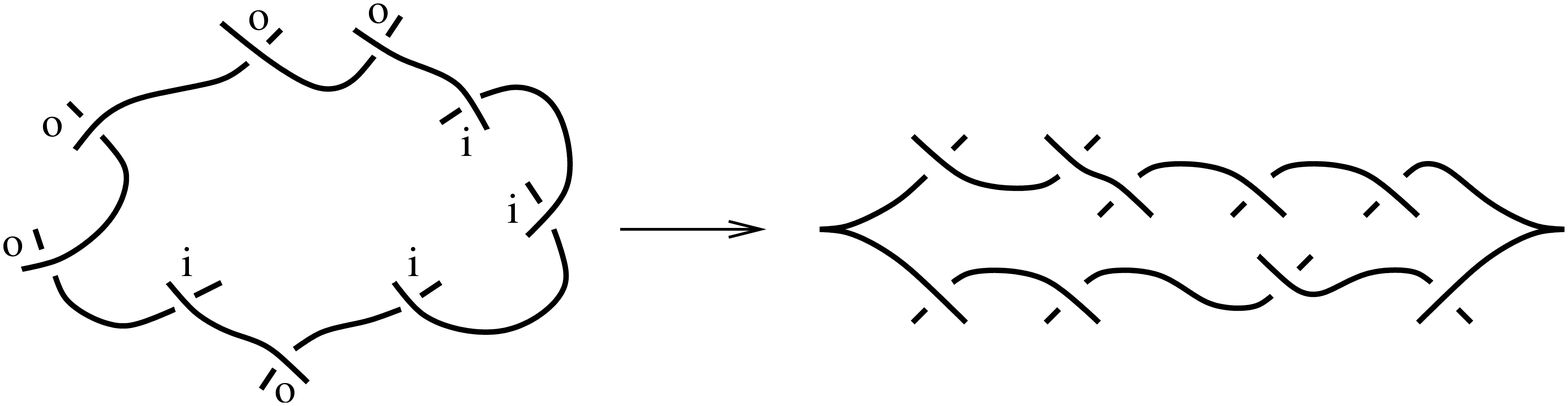}
   %\end{minipage}
   %\hfill
   %\begin{minipage}[c]{.5\textwidth}
   %\centering
   \caption{Outer (o) and inner (i) crossings along a state circle, and a hypothetical classification of the inner ones into three upper ones and one lower.}
   \label{fig:inandout}
   %\end{minipage}
\end{figure}

So in order to make the above idea precise, we need a slight strengthening of Proposition 11 of \cite{khov}. This is stated as Lemma \ref{lem:vazze}, and is proved after the current proof.

With it, once the first step (Mondrian diagram for the graph of the outermost state circles) is in place, all further steps in our process are carried out as follows. If a horizontal line segment corresponds to a state circle $s$ that contains other state circles, the outermost of which form the graph $G$, thicken the line segment into a rectangle $R$. Recall that the crossings of $D$ along $s$ that are inner to $s$ are classified into upper and lower ones. Augment $G$ with the vertices $u$ and $l$. If the state circle corresponding to the vertex $x$ of $G$ is connected to $s$ through an upper crossing, connect $x$ to $u$ with an edge, and similarly to $l$ if a lower crossing exists between the circle and $s$. To the resulting planar graph $G'$ with its two marked vertices, apply the construction of Lemma \ref{lem:vazze} to obtain the Mondrian diagram $M$ with upper- and lowermost horizontal line segments $u'$ and $l'$. After an isotopy, M can be copied into $R$ so that $u'$ and $v'$ (suitably lengthened) are the horizontal sides of $R$. The inner and outer crossings along $s=\partial R$ can be properly arranged by the following observation: Let $K$ be one of the half-infinite vertical strips above or below $R$. Any isotopy of $K$ to itself that fixes $\partial K$ pointwise and preserves horizontal and vertical lines takes iterated Mondrian diagrams to iterated Mondrian diagrams.

The planar diagram that results from the above inductive procedure can easily be turned into a front diagram by applying the operations shown in Figure \ref{fig:visszacsinal}: any horizontal line segments that weren't thickened into rectangles are now thickened into Legendrian unknot diagrams; rectangles are also isotoped into such; and finally, the vertical segments are changed into crossings.

\begin{figure}
   \centering
   \begin{minipage}[c]{.5\textwidth}
   \centering
   \includegraphics[width=\textwidth]{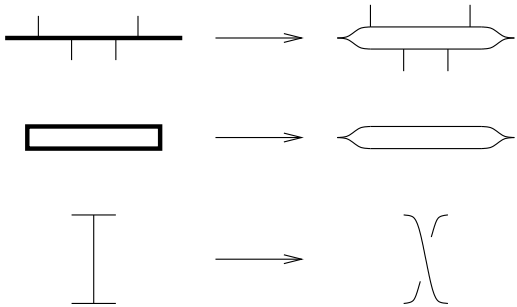}
   \end{minipage}
   \hfill
   \begin{minipage}[c]{.4\textwidth}
   \centering
   \caption{Turning an iterated Mondrian diagram into a front diagram}
   \label{fig:visszacsinal}
   \end{minipage}
\end{figure}

It's very easy to see that the front we constructed has an ungraded ruling: The unknot diagrams we've just constructed (Figure \ref{fig:visszacsinal}) from the state circles play the role of eyes and all crossings are switches. The corresponding discs are either disjoint or contained in one another, thus the normality condition is satisfied.
\end{proof}

\begin{lemma}\label{lem:vazze}
Let $G$ be a planar graph (possibly with multiple edges but without loop edges) with two marked vertices, $u$ and $l$, along the unbounded region of the complement of $G$. Then $G$ is the contraction of a Mondrian diagram so that the horizontal line segment $u'$ contracting to $u$ is the topmost one and the segment $l'$ contracting to $l$ is bottommost.
\end{lemma}

\begin{proof}
We will take advantage of some of the flexibilities left in Ng's construction. He first builds a `step-shaped' Mondrian diagram (see Figure 4 of \cite{khov}) for each so-called enhanced cycle $C$ of $G$. (An enhanced cycle is a subgraph of $G$ obtained as follows: Separate $G$ at its separating crossings. The (outer) boundary cycle of each resulting component, along with any other edges that connect two of its vertices, forms an enhanced cycle.) Instead of steps, we'll construct podium-shaped diagrams, as shown in Figure \ref{fig:lepcso}. The vertices that correspond to top and bottom (base) can be arbitrarily selected. The vertical pieces correspond to the edges of the boundary cycle $\tilde C$ of $C$. We choose the horizontal levels on the two sides so that they are all different. 

Now all edges in $C\setminus\tilde C$ have an upper and lower horizontal piece to connect. From each upper piece, we drop vertical line segments to the corresponding lower levels, arranged from left to right as suggested by the embedding of $C$ into the plane. Then we extend each left step to the right and each right step to the left to meet all the lower ends designated to it. The planarity of $C$ guarantees that this results in no `undesignated' intersections between horizontal and vertical segments, hence we get a Mondrian diagram. It is also strong in the sense of \cite{khov}, thus step 4 of Ng's construction does not have to be modified at all.

\begin{figure}
   \centering
   \begin{minipage}[c]{.6\textwidth}
   \centering
   \includegraphics[width=\textwidth]{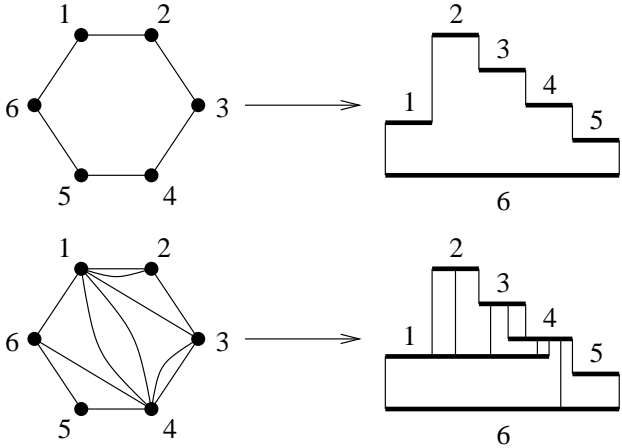} 
   \end{minipage}
   \hfill
   \begin{minipage}[c]{.3\textwidth}
   \centering
   \caption{Podium-shaped Mondrian diagrams for a cycle and for an enhanced cycle.}
   \label{fig:lepcso}
   \end{minipage}
\end{figure}

Then in step 3, Ng joins these building blocks in a tree-like fashion (Figure 5 of \cite{khov}), placing the base of each new step/podium onto the extension of the appropriate horizontal line segment. Let us use the enhanced cycle containing $l$ as the root of the tree and, of course, draw its podium-shaped diagram so that its base $l'$ contracts to $l$. 
%We also make sure that $u$ is not selected as a base for any of the enhanced cycles that contain it. 
As they are attached, draw the consecutive podium-shaped diagrams successively smaller so that at least one end of each horizontal piece is visible from above. (For segments that serve as a base at some point, this is achieved by making them stick out slightly from underneath the podium.) When we reach the enhanced cycle containing $u$ in the construction, we select the corresponding horizontal line segment $u'$ to be uppermost in its podium-shaped diagram. By the above visibility condition, $u'$ can be arranged to be the uppermost of all horizontal pieces thus far. In order to keep it that way, we need one last modification to Ng's construction: If at a later time, a podium-shaped diagram is to be attached to $u'$, we attach it from underneath, i.e.\ we construct the branch of the tree starting there upside down.
\end{proof}

\begin{kov}\label{szamolas}
In any $+$adequate diagram $D$ of the knot $A$, the writhe minus the number of components of $s_+(D)$ equals the minimum $v$--degree of the Kauffman polynomial minus one, which in turn equals $\TB(A)$. In the Dubrovnik version of the Kauffman polynomial, all coefficients of maximum $a$--degree terms are non-negative. Furthermore, these quantities agree with $\min\{\,k\,\big|\,\bigoplus_{i-j=k}HKh^{i,j}(K)\ne0\,\}$, where $HKh^{i,j}(K)$ is the Khovanov homology group in bigrading $(i,j)$.
\end{kov}

\begin{proof}
The front we constructed in Theorem \ref{front} has Thurston-Bennequin number as stated, which in turn is maximal and realizes the Kauffman bound by Rutherford's theorem. By the same, the said coefficients of the Kauffman polynomial represent counts of rulings and hence are not negative. The statement on Khovanov homology is a direct consequence of Proposition 8 of \cite{khov}.
\end{proof}

Thistlethwaite \cite{thist} observed that up to $11$ crossings, every knot or its mirror is $+$adequate.
%(???, see Stong). 
The invariant $\TB$ does distinguish mirrors, but at least for half of the knots up to $11$ crosings, Corollary \ref{szamolas} determines $\TB$.

\begin{pelda}\label{ex:KT}
The Kinoshita--Terasaka knot ($11n_{42}$) is a so-called adequate knot, i.e.\ it has a $+$adequate diagram which remains $+$adequate after switching all crossings. (Its trivial Alexander polynomial shows that the knot is not homogeneous \cite[section 7.6]{crom}, in particular it is neither alternating nor positive (or negative).) Thus, we may use our method to construct a maximum $tb$ diagram for both the knot and its mirror, showing that $\TB(11n_{42})=-7$ and $\TB(11n^*_{42})=-4$. (We obtain the same values for the Conway mutant of the knot ($11n_{34}$) and its mirror.)

\begin{figure}
   \centering
   \includegraphics[width=\linewidth]{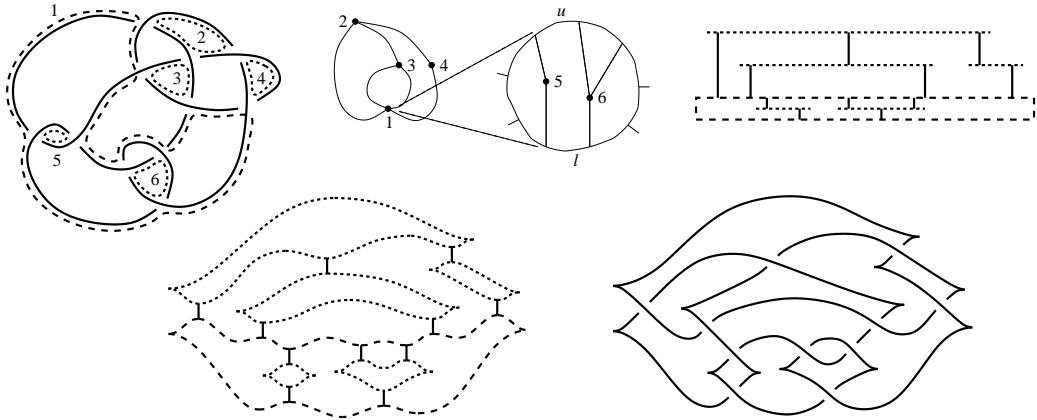} 
   \caption{Constructing a front diagram of the Kinoshita--Terasaka knot with maximum $tb=-7$.}
   \label{fig:kinoshita}
\end{figure}

\end{pelda}

\begin{figure}
   \centering
   \includegraphics[width=\linewidth]{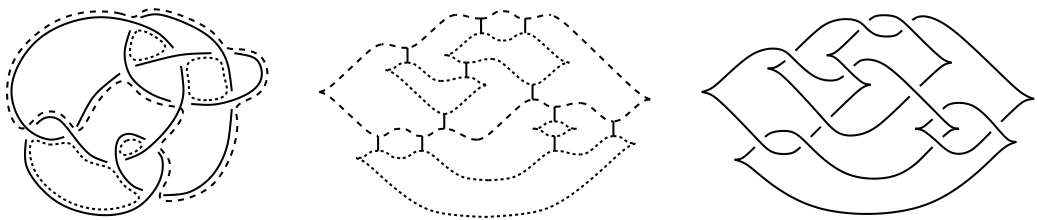} 
   \caption{Constructing a front diagram of the mirrored Kinoshita--Terasaka knot with maximum $tb=-4$.}
   \label{fig:terasaka}
\end{figure}

Acknowledgements: I would like to thank Lenny Ng, Tobias Ekholm, and Ko Honda for useful discussions.

\end{document}